\documentclass[reqno]{amsart}

\usepackage{amsmath,amssymb,amsthm}

\usepackage{tikz}

\usepackage{caption}
\usepackage{graphicx}
\usepackage{graphics}

\newtheorem*{thm}{Theorem}

\title[]{Kaczmarz Kac Walk}

\author[]{Stefan Steinerberger}
\address{Department of Mathematics, University of Washington, Seattle, WA 98195, USA} \email{steinerb@uw.edu}


\begin{document}

\begin{abstract} The Kaczmarz method is a way to iteratively solve a linear system of equations $Ax = b$. One interprets the solution $x$ as the point where hyperplanes intersect and then iteratively projects an approximate solution onto these hyperplanes to get better and better approximations. We note a somewhat related idea: one could take two random hyperplanes and project one into the orthogonal complement of the other. This leads to a sequence of linear systems $A^{(k)} x = b^{(k)}$ which is fast to compute, preserves the original solution and whose small singular values grow like
$\sigma_{\ell}(A^{(k)}) \sim \exp(k/n^2) \cdot \sigma_{\ell}(A)$.
\end{abstract}

 \maketitle

\section{Introduction and Result}
\subsection{The Kaczmarz method.} The Kaczmarz method \cite{kac} is an iterative method for solving linear systems of equations proposed in 1937. Let $A \in \mathbb{R}^{m \times n}$, $m \geq n$, be a given matrix, either square or overdetermined. We assume that $A:\mathbb{R}^n \rightarrow \mathbb{R}^m$ is injective and that
$ Ax =b,$ where $x \in \mathbb{R}^n$ is unknown and $b \in \mathbb{R}^m$ is a known right-hand side. Using $A_1, \dots, A_m \in \mathbb{R}^n$ to denote the rows of $A$, Kaczmarz proposes to interpret the linear system $Ax = b$ geometrically and to think of 
the solution $x$ as the unique point in the intersection of the hyperplanes $H_1, \dots, H_m \subset \mathbb{R}^m$ where
$$ H_i = \left\{w \in \mathbb{R}^n: \left\langle A_i, w \right\rangle = b_i \right\}  \qquad \mbox{for}~1 \leq i \leq m.$$

\vspace{-10pt}

\begin{center}
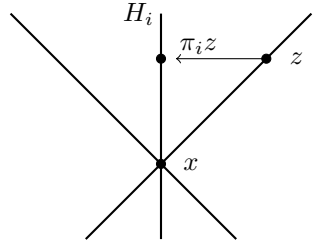
\begin{figure}[h!]
\begin{tikzpicture}[scale=2]
\draw[thick] (-0.5, -0.5) -- (1, 1);
\draw[thick] (0, -0.5) -- (0, 1);
\draw[thick] (0.5, -0.5) -- (-1, 1);
\filldraw (0,0) circle (0.03cm);
\node at (0.2, 0) {$x$};
\filldraw (0.7, 0.7) circle (0.03cm);
\node at (0.9, 0.7) {$z$};
\filldraw (0, 0.7) circle (0.03cm);
\node at (0.25, 0.78) {$\pi_i z$};
\draw [ <-] (0.1, 0.7) -- (0.7, 0.7);
\node at (-0.15, 1) {$H_i$};
\end{tikzpicture}
\caption{Projection $\pi_i y$ onto $H_i$ given by $\left\langle a_i, w\right\rangle = b_i$.}
\label{fig:1}
\end{figure}
\end{center}
\vspace{-20pt}
The Pythagorean Theorem implies that, if $z \in \mathbb{R}^n$ is any point and $\pi_i:\mathbb{R}^n \rightarrow H_i$ is the orthogonal projection onto any hyperplane, then $ \|x - \pi_i z \| \leq \|x-z\|$
with equality if and only if $z \in H_i$ (see Fig. \ref{fig:1}). By iteratively applying projections onto the $m$ hyperplanes, we iteratively decrease the distance to the solution $x$. Each individual projection is fast to compute since
$$\pi_i y = y+ \frac{b_i - \left\langle A_i, y\right\rangle}{\|A_i\|^2}A_i.$$
This has made the method popular when dealing with very large matrices: it is a row-based approach that does not require to load the entire matrix into memory. The main remaining question is in which order one should choose the projections/rows. Here, an elegant solution was given by Strohmer--Vershynin \cite{strohmer} who proposed to pick the projections  randomly with likelihood proportional to the size of the row when interpreted as a vector in $\ell^2(\mathbb{R}^n)$. 

\begin{thm}[Strohmer--Vershynin] If the projection onto hyperplane $H_i$ is chosen with likelihood $\|A_i\|^2/\|A\|_F^2$, then
\begin{align*} 
\mathbb{E}~ \|x_k - x\|^2 \leq \left(1 - \frac{\sigma_{\min}(A)^2}{\|A\|_F^2}\right)^k \|x_0 - x\|^2,
\end{align*}
where $\sigma_{\min}$ is the smallest singular value of $A$.
\end{thm}
 This rate is known to be optimal \cite{stein}. The Kaczmarz method has been used in a wide variety of settings, we mention early related work of Agmon \cite{agmon}, Motzkin \cite{motzkin} and more recent work  \cite{alder, bai0, bai, gower, hadd, kac,leventhal, man, ma, mar, need2, stein2, stein3, stein4, stein5, stein6, zouz}.
  
 \subsection{Kaczmarz Kac Walk}
A basic idea when trying to solve a linear system of equations $Ax = b$ is to replace it by another system $A^* x = b^*$ which has the same solution but a `better' matrix $A^*$. Returning to the main idea behind the Kaczmarz method, suppose we have two hyperplanes
$$ \left\langle A_i,x \right\rangle = b_i \qquad \mbox{and} \qquad  \left\langle A_j,x \right\rangle = b_j$$
with rows normalized to length 1, i.e. $\|A_i\| = 1 = \|A_j\|$. One could take a linear combination of the two and replace any of the two equations by that linear combination. A particularly natural choice is to project one of the hyperplanes onto the orthogonal complement of the other and then renormalizing, that is
\begin{equation} \label{new}
 \left\langle \frac{A_j - \left\langle A_i, A_j \right\rangle  A_i }{\| A_j - \left\langle A_i, A_j \right\rangle  A_i \|},x \right\rangle = \frac{b_j - \left\langle A_i, A_j \right\rangle b_i}{\|A_j - \left\langle A_i, A_j \right\rangle  A_i\|}.
 \end{equation}
 This obviously requires $A_i$ and $A_j$
to be linearly independent or, since both are normalized to $\|A_i\| = 1 = \|A_j\|$, we only need $A_i \neq \pm A_j$.
Equation \eqref{new} requires the computation of a single inner product: $A_i$ and $ A_j - \left\langle A_i, A_j \right\rangle  A_i$ are orthogonal, the Pythagorean Theorem implies 
$$\| A_j - \left\langle A_i, A_j \right\rangle  A_i \| = \sqrt{1 - \left\langle A_i, A_j \right\rangle^2}.$$
This suggests the following procedure.
\begin{quote}
\textbf{Kaczmarz Kac Walk.} For $A \in \mathbb{R}^{m \times n}$ with rows normalized to length $1$, pick two rows $i \neq j$ uniformly at random and if $A_i \neq \pm A_j$,
\begin{equation} \label{main}
 A_j \leftarrow \frac{ A_j - \left\langle A_i, A_j \right\rangle A_i}{\sqrt{1 - \left\langle A_i, A_j \right\rangle^2}}
 \end{equation} 
 as well as $b_j \leftarrow (b_j - \left\langle A_i, A_j \right\rangle b_i)/ \sqrt{1 - \left\langle A_i, A_j \right\rangle^2}$.
\end{quote}
This process should make the rows `more orthogonal' and the new linear system should be better conditioned and easier to solve.  Some quick numerical experiments show that this indeed happens (see Fig. 1). The question is now whether and in what sense the improvement can be quantified.

 \subsection{Main Result}
We now discuss the main result: the procedure is beneficial insofar as it regularizes the matrix along degenerate subspaces where the singular values are $<1$ (in expectation) while keeping the Frobenius norm invariant.

\begin{thm}
Let $A \in \mathbb{R}^{m \times n}$ be invertible with normalized rows, $\|A_i\|_{\ell^2}  = 1$, that are different ($A_i \neq \pm A_j$ when $i \neq j$). Use $\phi(A) \in \mathbb{R}^{m \times n}$ to denote the matrix obtained by selecting $i \neq j$ uniformly at random from $\left\{1,2,\dots,m\right\}$ and applying \eqref{main}.  Then,
$$ \forall x \in \mathbb{R}^n \qquad    \mathbb{E}~ \|\phi(A) x\|^2 \geq \|Ax\|^2 +  \frac{2}{m(m-1)} \left( \|Ax\|^2 - \|A^T Ax\|^2 \right).$$
\end{thm}
We first discuss the case of square matrices, $m=n$, that setting is slightly simpler. 
The procedure maintains the normalization $\|A_i\|_{\ell^2} = 1$ and thus, in particular, the squared sum of the singular values remains invariant
$$ \sum_{i=1}^{n} \sigma_i(A)^2 = n =  \sum_{i=1}^{n} \sigma_i(\phi(A))^2.$$
This invariance implies that $\phi(A)$ is not going to be uniformly larger than $A$, we should expect $\mathbb{E} \|\phi(A) x\|^2$ to be sometimes larger and sometimes smaller than $\|A x\|^2$, it should depend on $x$. This is exactly what happens.

\begin{center}
\begin{figure}[h!]
\begin{tikzpicture}
\node at (0,0) {\includegraphics[width=0.45\textwidth]{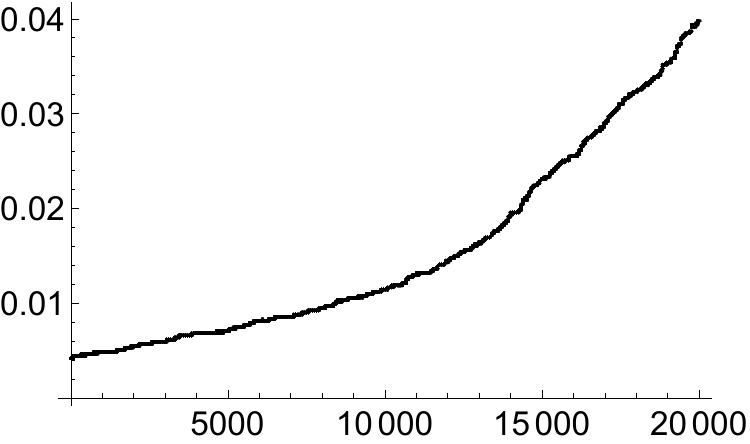}};
\node at (6,0) {\includegraphics[width=0.45\textwidth]{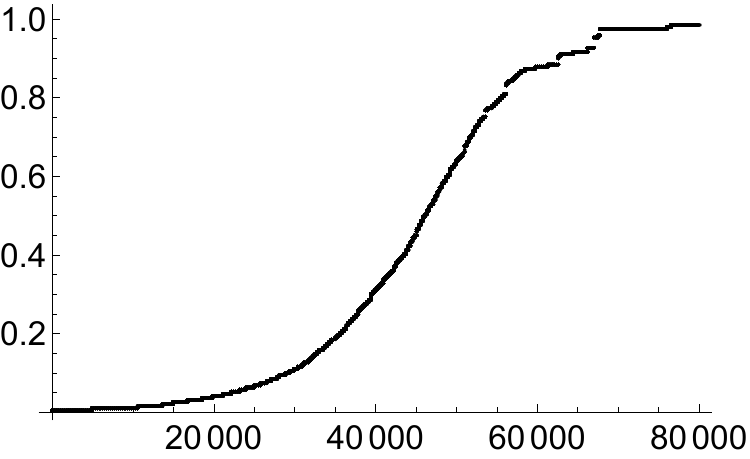}};
\end{tikzpicture}
\caption{The smallest singular value of random Gaussian $100 \times 100$ matrix over 20000 (left) and 80000 (right) iterations.}
\label{fig:2}
\end{figure}
\end{center}
\vspace{-10pt}
Let $x = v_n$ be a singular vector associated to the smallest singular value $\sigma_n$ of $A$. The spectral normalization implies $\sigma_n \leq 1$ with equality if and only if $A$ is an orthogonal matrix. Then, if $\sigma_n < 1$,
$$ \mathbb{E}~\| \phi(A) v_n\|^2 \geq \sigma_n^2 + \frac{2}{n(n-1)}\left( \sigma_n^2 - \sigma_n^4\right) > \sigma_n^2.$$
The observation can be generalized: let us consider all singular values
$ \sigma_1 \geq \sigma_2 \geq \dots \sigma_\ell \geq 1 > \sigma_{\ell+1} \geq \dots \geq \sigma_n$
with the associated singular vectors $v_1, v_2, \dots, v_n$. If $x \in \mbox{span}\left\{v_{\ell+1}, \dots, v_n\right\}$ is in the `bad' (difficult to invert) subspace, then
$$ \|Ax\|^2 = \left\| A \sum_{k=\ell+1}^{n} a_k v_k \right\|^2 = \sum_{k= \ell+1}^{n} a_k^2 \sigma_k^2 \geq  \sum_{k= \ell+1}^{n} a_k^2 \sigma_k^4 = \|A^T Ax\|^2.$$

\begin{center}
\begin{figure}[h!]
\begin{tikzpicture}
\node at (0,0) {\includegraphics[width=0.45\textwidth]{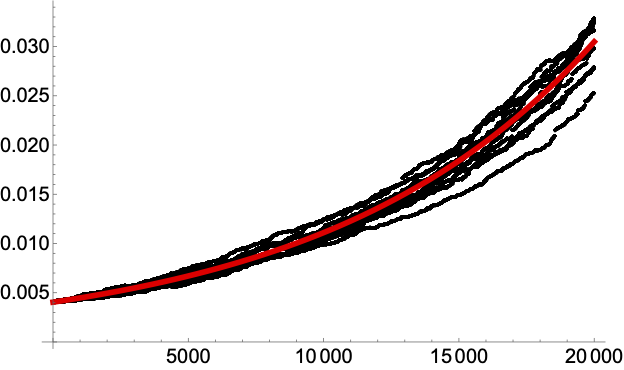}};
\node at (6,0) {\includegraphics[width=0.45\textwidth]{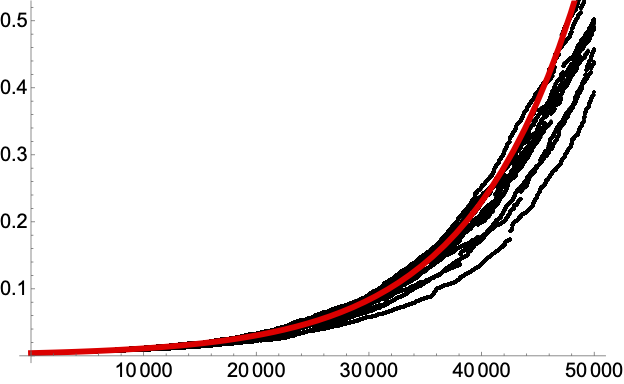}};
\end{tikzpicture}
\caption{Prediction \eqref{pred} in red tested against ten random evolutions of $\sigma_{100}(A^{(k)})$
for a random initial matrix of size $100 \times 100$. We observe high accuracy while
$\sigma_{100}(A^{(k)}) \ll 1$.}
\label{fig:3}
\end{figure}
\end{center}

This means that $\phi(A)$ expands relative to $A$ along directions where $A$ is `small', meaning the subspace spanned by singular vectors corresponding to singular values $<1$. The same argument, with $\geq$ replaced by $\leq$, shows that $\phi(A)$ contracts relative to $A$ along directions where $A$ is `large', that is the subspaces spanned by singular vectors corresponding to singular values $>1$.  One would expect behavior like this to eventually lead to an orthogonal matrix and this is what is being observed empirically, see Fig. \ref{fig:2}.
  More can be said. If $\sigma_{\ell}(A) \ll 1$, then the result suggests the approximation
$$ \sigma_{\ell}(\phi(A))^2 \sim \mathbb{E}~\| \phi(A) v_{\ell}\|^2 \geq \left(1 + \frac{2}{n(n-1)}\right)\sigma_{\ell}(A)^2$$
which, in turn, suggests that for the sequence of matrices $A = A^{(0)}, A^{(1)}, A^{(2)}, \dots$ arising from the Kaczmarz preconditioner, one would expect
\begin{align} \label{pred}
 \sigma_{\ell}(A^{(k)}) \sim \left(1 + \frac{2}{n(n-1)}\right)^{k/2}\sigma_{\ell}(A) \sim  e^{k/n^2} \cdot \sigma_{\ell}(A). 
 \end{align} 
It is clear that one can only hope for this simple prediction to be accurate as long as $\sigma_{\ell} \ll 1$ since we are ignoring other terms. In that regime, however, the prediction is quite accurate (see Fig. \ref{fig:3}).

\begin{center}
\begin{figure}[h!]
\begin{tikzpicture}
\node at (0,0) {\includegraphics[width=0.45\textwidth]{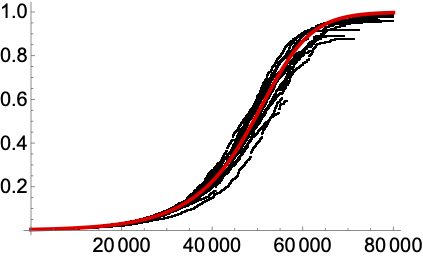}};
\node at (6,0) {\includegraphics[width=0.45\textwidth]{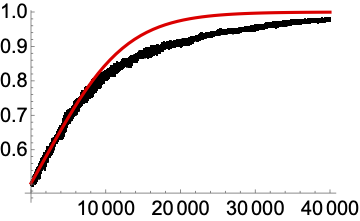}};
\end{tikzpicture}
\caption{Prediction \eqref{pred3} tested against ten random evolutions of $\sigma_{100}(A^{(k)})$ (left) and $\sigma_{70}(A^{(k)}$ (right) for a random $100 \times 100$ matrix.}
\label{fig:34}
\end{figure}
\end{center}

More generally, we can turn the heuristic
$$ \sigma_{\ell}(\phi(A))^2 \sim \mathbb{E}~\| \phi(A) v_{\ell}\|^2 \geq \sigma_{\ell}^2 + \frac{2}{n(n-1)}\left( \sigma_{\ell}^2 - \sigma_{\ell}^4\right)$$
into an ordinary differential equation for $y(t) = \sigma_{\ell}(A^{\left\lfloor t \right\rfloor})^2$ which reads
$$ y'(t) = \frac{2}{n(n-1)} (y(t) - y(t)^2) \quad \mbox{and solved by} \quad y(t) = \left(1 + c \cdot e^{-\frac{2t}{n(n-1)}}\right)^{-1}.$$
Solving for the initial condition $y(0) = \sigma_{\ell}(A)^2$, we arrive at the prediction
\begin{align} \label{pred3}
 \sigma_{\ell}(A^{(k)}) \sim \left(1 + \left( \frac{1}{\sigma_{\ell}(A)^2} - 1 \right) \cdot e^{-\frac{2k}{n(n-1)}}\right)^{-1/2}
 \end{align} 
 This prediction accounts for the nonlinear ODE and remains accurate for much longer: examples are shown in Fig. \ref{fig:34}. Prediction \eqref{pred3} is reasonably accurate for the smallest singular values and less accurate for larger ones; one possible interpretation is that the procedure introduces rotation of vectors associated to smaller singular values which allows for a form of `leakage' to other singular vectors; it would be interesting to have a better understanding of this.
 
 \subsection{The overdetermined case.} 
\begin{center}
\begin{figure}[h!]
\begin{tikzpicture}
\node at (0,0) {\includegraphics[width=0.45\textwidth]{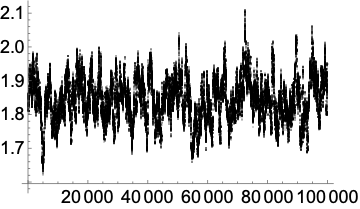}};
\node at (6,0) {\includegraphics[width=0.35\textwidth]{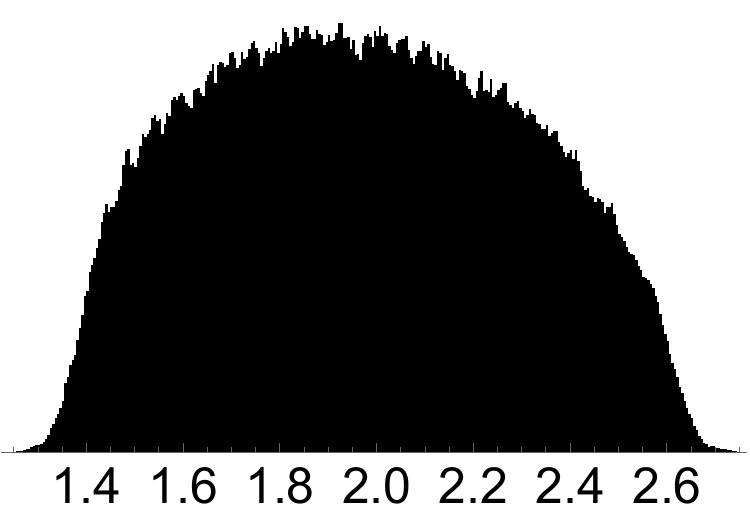}};
\end{tikzpicture}
\caption{Evolution of a $100 \times 25$ matrix: $\sigma_1(A^{(k)})/\sigma_{25}(A^{(k)})$ (left) and empirical density of all singular values (right).}
\label{fig:4}
\end{figure}
\end{center}

The over-determined case $m > n$ is more complicated. The Theorem still applies as stated, however, the arising behavior is more complicated. A first indication of this is that the spectral invariant now reads
 $$ \sum_{i=1}^{n} \sigma_i(A)^2 = m =  \sum_{i=1}^{n} \sigma_i(\phi(A))^2$$
which rules out the possibility of all singular values being 1. This is not surprising, we have $m>n$ vectors in $\mathbb{R}^n$, some will not be mutually orthogonal. The Theorem applies and suggests growth of small singular values in expectation while the overall Frobenius norm, the sum of the squared singular values, stays invariant: small singular values expand, large singular values contract. An example is shown in Fig. \ref{fig:4} which shows the evolution of an overdetermined $100 \times 25$ matrix. The condition number remains reasonable small; moreover, the singular values seem to follow an underlying density. It would be interesting to see whether this distribution can be identified and whether the convergence of the empirical singular values to this limiting distribution can be rigorously established.

\begin{center}
\begin{figure}[h!]
\begin{tikzpicture}
\node at (0,0) {\includegraphics[width=0.95\textwidth]{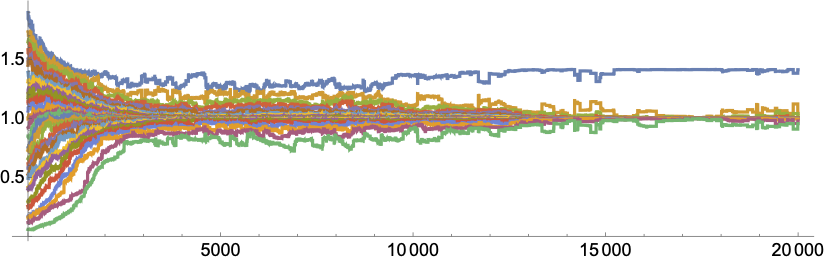}};
\end{tikzpicture}
\caption{Evolution of a $31 \times 30$ matrix. The picture suggests that there is a tendency towards $\sigma_1 \sim \sqrt{2} \geq 1 \sim \sigma_2 \sim \dots \sigma_n$.}
\label{fig:5}
\end{figure}
\end{center}

There is a curious special case that deserves to be emphasized, the case $m=n+1$. If we take a matrix $A \in \mathbb{R}^{3 \times 2}$ with rank 2, then the procedure will eventually return an orthonormal set of 2 vectors with a third vector of size 1, the singular values will be $\left\{\sqrt{2},1\right\}$. Something similar seems to happen for matrices in $A \in \mathbb{R}^{(n+1) \times n}$, the procedure seems to shift them to an orthonormal basis with one extra vector and thus a matrix with singular values $\left\{\sqrt{2}, 1, 1, \dots, 1\right\}$. An example is shown in Fig. \ref{fig:5}. This is quite interesting, we do not know why this happens.

\vspace{-20pt}
The behavior of matrices in $\mathbb{R}^{(n+2) \times n}$ appears to be fundamentally different, no such limiting behavior seems to emerge, or maybe it simply requires a much longer timescale to emerge. The behavior of matrices in $\mathbb{R}^{n \times 2}$ is also somewhat amusing (see also \S 2.3): since each row is normalized in $\ell^2$, we can identify these matrices with a collection of $n$ points on $\mathbb{S}^1$. The procedure seems to lead to a clustering of the points along the four directions $\pm e_1, e_2$ or a rotation thereof. Numerical issues, two vectors being numerically indistinguishable, arise quickly. In contrast, the case $\mathbb{R}^{n \times 3}$, identifiable with $n$ points on $\mathbb{S}^2$, appears to be perfectly random throughout.

\subsection{Kac walk.} We are not aware of this procedure having been previously investigated, however, it is obviously similar Kac's random walk introduced in 1956 \cite{kac} as a toy model for the Boltzmann equation. Kac assumes we are given $n$ points on $\mathbb{S}^{n-1}$, identifiable with a matrix $A \in \mathbb{R}^{n \times n}$ whose rows are normalized to have $\ell^2-$norm 1. At the $k-$th step, one then picks two different rows $1 \leq i \neq j \leq n$ uniformly at random and updates, for $0 \leq \theta \leq 2\pi$ chosen uniformly at random,
$$ \begin{pmatrix} A_i(k+1) \\ A_j(k+1) \end{pmatrix} = \begin{pmatrix} \cos{(\theta)} & \sin{(\theta)} \\ - \sin{(\theta)} & \cos{(\theta)} \end{pmatrix} \begin{pmatrix} A_i(k) \\ A_j(k) \end{pmatrix}.$$
This corresponds to a random rotation in the two-dimensional subspace generated by the rows $A_i$ and $A_j$. The idea is that this process should suitable randomize system and one should get the uniform distribution; this has now been well understood \cite{jain, jiang, mischler, pillai} at the quantitative level. It is clear that our model is substantially different: for square matrices, for example, we do not expect uniform distribution but convergence to an orthogonal matrix. A remaining question is when the Kaczmarz Kac Walk can be useful from a numerical perspective. It really does replace the system by a system whose condition number improves, a system that is easier to solve. Moreover, it effectively improves all small singular values simultaneously. It is clear that if we have a sufficiently large computational budget, it will always be worth it to spend some of it on the Kaczmarz Kac walk (as one then inherits the improved condition number later on); whether this, or some natural variation thereof, may turn out to be numerically useful remains to be seen.

 \section{Proof of the Theorem}
 \subsection{Proof of the Theorem}
 \begin{proof}
Let us fix some $x \in \mathbb{R}^n$. Suppose we choose $i \neq j$ uniformly at random and replace the $i-$th equation by the $i-$th equation after having been projected onto the orthogonal complement of the $j-$th equation and then rescaled, formally
$$ A_i \leftarrow   \frac{A_i - \left\langle A_i, A_j \right\rangle A_j}{\sqrt{1 - \left\langle A_i, A_j\right\rangle^2}}.$$
Let us denote the new matrix by $A_2$. Our goal is to understand how $\|A_2x\|^2$ changes compared to $\|Ax\|^2$. 
Since we only modify the $i-$th row of the matrix, the change is completely explicit and we have
$$ (A_2 x)_k = \begin{cases} \frac{(Ax)_i - \left\langle A_i, A_j \right\rangle (Ax)_j }{\sqrt{1 - \left\langle A_i, A_j\right\rangle^2}} \qquad &\mbox{if}~k=i \\
(Ax)_k \qquad &\mbox{otherwise.} \end{cases}$$
Therefore, summing over all possible combinations of $1 \leq i,j \leq m$ with $i \neq j$,
\begin{align*}
\mathbb{E} \|A_2 x\|^2 &= \frac{1}{m(m-1)}\sum_{i,j=1 \atop i \neq j}^m \|A_2 x\|^2 \\
&=   \frac{1}{m(m-1)}\sum_{i,j=1 \atop i \neq j}^m \|Ax\|^2 - \left\langle A_i, x\right\rangle^2+ \left( \frac{ \left\langle A_i, x\right\rangle - \left\langle A_i, A_j \right\rangle \left\langle A_j, x\right\rangle }{\sqrt{1 - \left\langle A_i, A_j\right\rangle^2}} \right)^2.
\end{align*}
We note that $\|Ax\|$ does not depend on $i$ or $j$, one has
 $$ \mathbb{E} \|A_2 x\|^2 = \|Ax\|^2 + \frac{1}{m(m-1)}\sum_{i,j=1 \atop i \neq j}^m  -  \left\langle A_i, x\right\rangle^2+ \left( \frac{ \left\langle A_i, x\right\rangle - \left\langle A_i, A_j \right\rangle  \left\langle A_j, x\right\rangle }{\sqrt{1 - \left\langle A_i, A_j\right\rangle^2}} \right)^2.$$
It remains to control the sum 
$$ \Sigma = \sum_{i,j=1 \atop i \neq j}^m  -  \left\langle A_i, x\right\rangle^2+ \left( \frac{ \left\langle A_i, x\right\rangle - \left\langle A_i, A_j \right\rangle  \left\langle A_j, x\right\rangle }{\sqrt{1 - \left\langle A_i, A_j\right\rangle^2}} \right)^2$$
and the remainder of the proof is dedicated to that. We observe that, trivially,
$$ \frac{1}{1 - \left\langle A_i, A_j\right\rangle^2} \geq 1$$
which is attained when $m=n$ and the matrix $A$ is orthogonal. It is clear that when $m > n$, this step in the argument is going to be lose nontrivial factors, it would be nice if these could be recovered, perhaps by means of a different argument. Using this trivial bound, we obtain the lower bound
 $$ \Sigma \geq \sum_{i,j=1 \atop i \neq j}^m  -  \left\langle A_i, x\right\rangle^2+  \left( \left\langle A_i, x\right\rangle - \left\langle A_i, A_j \right\rangle  \left\langle A_j, x\right\rangle \right)^2.$$
Expanding the square and simplifying, we arrive at
  \begin{align*}
  \Sigma &\geq  \sum_{i,j=1 \atop i \neq j}^m  -  \left\langle A_i, x\right\rangle^2+  \left( \left\langle A_i, x\right\rangle - \left\langle A_i, A_j \right\rangle  \left\langle A_j, x\right\rangle \right)^2 \\
  &=  \sum_{i,j=1 \atop i \neq j}^m  - 2 \left\langle A_i, x\right\rangle\left\langle A_i, A_j \right\rangle  \left\langle A_j, x\right\rangle  + \left\langle A_i, A_j \right\rangle^2  \left\langle A_j, x\right\rangle^2.
  \end{align*}
  It will be easier to work with this expression if we remove the restriction $i \neq j$. When $i =j$, the summand simplifies to
  $$ - 2 \left\langle A_i, x\right\rangle\left\langle A_i, A_i \right\rangle  \left\langle A_i, x\right\rangle  + \left\langle A_i, A_i \right\rangle^2  \left\langle A_i, x\right\rangle^2 = -\left\langle A_i, x\right\rangle^2$$
  and thus adding over the diagonal leads to an additional $-\|Ax\|^2$ term. Hence
  $$ \Sigma \geq \|Ax\|^2 + \sum_{i,j=1}^m  - 2 \left\langle A_i, x\right\rangle\left\langle A_i, A_j \right\rangle  \left\langle A_j, x\right\rangle  + \left\langle A_i, A_j \right\rangle^2  \left\langle A_j, x\right\rangle^2.$$
 We treat both summands separately. Abbreviating $y = Ax$, we may think of the first summand as a quadratic form and have
 \begin{align*}
 \sum_{i,j=1}^m \left\langle A_i, x\right\rangle\left\langle A_i, A_j \right\rangle  \left\langle A_j, x\right\rangle &= \sum_{i,j=1}^m y_i \left\langle A_i, A_j \right\rangle  y_j \\
 &= \left\langle y, A A^T y \right\rangle = \left\langle Ax, A A^T Ax \right\rangle \\
 &=  \left\langle A^T Ax, A^T Ax \right\rangle = \|A^T Ax\|^2
 \end{align*}
and the first summand can be written as
  $$ -2 \sum_{i,j=1}^m   \left\langle A_i, x\right\rangle\left\langle A_i, A_j \right\rangle  \left\langle A_j, x\right\rangle  = -2 \|A^TAx\|^2.$$
  The second summand is non-negative and by only summing over $i=j$, we get
 $$  \sum_{i,j=1}^m \left\langle A_i, A_j \right\rangle^2  \left\langle A_j, x\right\rangle^2 \geq
 \sum_{j=1}^{m} \left\langle A_j, x\right\rangle^2 = \|Ax\|^2.$$
Collecting all the terms, we have
$$ \frac{\Sigma}{m(m-1)}  \geq \frac{2}{m(m-1)} \left( \|Ax\|^2 -\|A^T Ax\|^2 \right).$$
  \end{proof}

  \subsection{Remark.} One could use, valid for all $-1 < x < 1$,
$$ \frac{1}{1-x^2} \geq 1+x^2$$
to argue that
\begin{align*}
 \left( \frac{ \left\langle A_i, x\right\rangle - \left\langle A_i, A_j \right\rangle  \left\langle A_j, x\right\rangle }{\sqrt{1 - \left\langle A_i, A_j\right\rangle^2}} \right)^2 &=  \frac{\left( \left\langle A_i, x\right\rangle - \left\langle A_i, A_j \right\rangle  \left\langle A_j, x\right\rangle\right)^2 }{1 - \left\langle A_i,A_j\right\rangle^2} \\
 &\geq \left( \left\langle A_i, x\right\rangle - \left\langle A_i, A_j \right\rangle  \left\langle A_j, x\right\rangle\right)^2 (1 + \left\langle A_i,A_j\right\rangle^2)
 \end{align*}
to derive the presence of an additional factor
    $$ \Sigma_2 = \sum_{i,j=1 \atop i \neq j}^m  \left( \left\langle A_i, x\right\rangle - \left\langle A_i, A_j \right\rangle  \left\langle A_j, x\right\rangle \right)^2  \left\langle A_i, A_j \right\rangle^2.$$
  Where as the main proof showed that
   $$ \mathbb{E} \|A_2 x\|^2 \geq \|Ax\|^2 +  \frac{\Sigma}{m(m-1)}$$
  and the bounded $\Sigma$ from below in terms of $ \|Ax\|$ and $\|A^T Ax\|$, this slightly refined argument would imply
     $$ \mathbb{E} \|A_2 x\|^2 \geq \|Ax\|^2 +  \frac{\Sigma + \Sigma_2}{m(m-1)}.$$
It is clear that $\Sigma_2$ is a `higher-order term' insofar as larger powers of $A$ play a role; this means that we can expect it to be negligible when $m=n$ (as also evidenced by the numerics, see above). However, in the over-determined regime $m>n$, this is no longer the case.

\subsection{The case $A \in \mathbb{R}^{n \times 2}$}
 We conclude with some observations for the case $A \in \mathbb{R}^{n \times 2}$. These observations are non-generic in the sense that this case appears to be slightly too rigid to be in any way representative of the general case. 
 What one observes numerically is that the $n$ points, which we can identify with $n$ points on $\mathbb{S}^1$ or, even more conveniently, by their angles $\theta \in [0,2\pi]$ end up in 4 clusters at the angles $x_0, x_0 + \pi/2, x_0 + \pi$ and $x_0 + 3 \pi/2$. This is interesting as it corresponds to some type of singularity formation. We have no proof of this observation.
 The remainder of the section is dedicating to proposing a mean-field limit. We identify the $n \gg 1$ points with their angles in $[0,2 \pi]$ and assume that $n$ is so large that it makes sense to think of them as described by a smooth probability density $u:[0,2\pi] \rightarrow \mathbb{R}_{\geq 0}$ of points.  The entire problem is invariant under translation. 
 Consider what happens to the density in $[-\varepsilon, \varepsilon]$ for some small $\varepsilon > 0$.
  If we pick a point with angle in $[-\varepsilon, \varepsilon]$, it is very likely to be transported somewhere else, so we expected 
  $$ \mbox{a constant loss proportional to} \qquad \int_{-\varepsilon}^{\varepsilon} u(x) dx.$$
  However, we may also pick another point, at angle $s$, and by suitable projection onto the orthogonal line corresponding to a second point at angle $t$, send it to the interval $[-\varepsilon, \varepsilon]$. Suppose we pick  $(\cos{s}, \sin{s})$ and project it onto the orthogonal complement of $(\cos{t}, \sin{t})$: then
$$ \frac{\begin{pmatrix} \cos{s} \\ \sin{s} \end{pmatrix} - \left\langle  \begin{pmatrix} \cos{s} \\ \sin{s} \end{pmatrix},  \begin{pmatrix} \cos{t} \\ \sin{t} \end{pmatrix} \right\rangle  \begin{pmatrix} \cos{t} \\ \sin{t} \end{pmatrix}}{\left\|\begin{pmatrix} \cos{s} \\ \sin{s} \end{pmatrix} - \left\langle  \begin{pmatrix} \cos{s} \\ \sin{s} \end{pmatrix},  \begin{pmatrix} \cos{t} \\ \sin{t} \end{pmatrix} \right\rangle  \begin{pmatrix} \cos{t} \\ \sin{t} \end{pmatrix}\right\|} = \frac{1}{|\sin{(s-t)}|} \begin{pmatrix} - \sin{(t)} \sin{(s-t)} \\ \cos{(t)} \sin{(s-t)} \end{pmatrix}$$
If we want this to correspond to an angle in $[-\varepsilon, \varepsilon]$, we need $|\cos{t}|$ to be small and thus (ignoring lower order terms in $\varepsilon^2$)
$$ t \in \left[ \frac{\pi}{2} - \varepsilon, \frac{\pi}{2} + \varepsilon\right] \cup \left[ \frac{3\pi}{2} - \varepsilon, \frac{3\pi}{2} + \varepsilon\right].$$
Having $t$ in that set ensures the point being sent to either $[-\varepsilon, \varepsilon]$ or $[\pi - \varepsilon, \pi+\varepsilon]$. To ensure it landing in $[-\varepsilon, \varepsilon]$, we need $-\pi/2 \leq s \leq \pi/2$ (if $t \sim \pi/2$) and $\pi/2 \leq s \leq 3\pi/2$ (if $t \sim 3\pi/2$) and it does not matter where in that interval $s$ lands.
Collecting these main contributions suggests that the mean-field limit might be governed by the integro-differential equation
\begin{align*}
 \frac{\partial }{\partial t} u(t,x) &= -u(t,x) \\
 &+ u(t, x + \pi/2) \int_{x - \pi/2}^{x+\pi/2} u(t,y) dy +  u(t, x - \pi/2) \int_{x + \pi/2}^{x+3\pi/2} u(t,y) dy.
 \end{align*}
The constant density $u(x) = (2\pi)^{-1}$ is a steady-state solution of this equation; however, suggested by the empirically observed singularity formation, it should not be stable under perturbations. Making the ansatz $ u(t,x) = 1/(2\pi) + \varepsilon f(t,x)$ for some small $\varepsilon$ and function $f(t,x)$ with mean value 0, one has, to leading order,
\begin{align*}
  \frac{\partial }{\partial t} f(t,x) &= -  f(t,x) + \frac{1}{2} f(t,x-\pi/2) + \frac{1}{2} f(t,x+\pi/2) 
\end{align*}  
Separation of variables shows that this non-local equation has the solutions
  $$ f(t,x) = \sum_{k \geq 1} e^{-2 \sin\left( k \pi/4 \right)^2 t} \left(a_k \sin{(kx)} + b_k \cos{(kx)}\right).$$
This is an unusual type of stability that dampens all frequencies except those that are a multiple of 4, in which case $\sin{(k\pi/4)}^2 = 0$. This is consistent with the formation of four clusters at four equispaced angles, however, it also indicates that the phenomenon may arise from a higher-order nonlinearity.

\end{document}